\documentclass{amsart}
\usepackage{color}

\usepackage{amsmath} 
\usepackage{amssymb}
\usepackage{mathrsfs}

\newtheorem{theorem}{Theorem}[section] 
\newtheorem{claim}[theorem]{Claim}

\newtheorem{conclusion}[theorem]{Conclusion}

\theoremstyle{definition}
\newtheorem{definition}[theorem]{Definition}

\newtheorem{problem}[theorem]{Problem}

\newtheorem{discussion}[theorem]{Discussion}

\theoremstyle{remark}
\newtheorem{remark}[theorem]{Remark}

\newtheorem{notation}[theorem]{Notation}
\newtheorem{context}[theorem]{Context}

\newcommand{\lf}{{\rm lf}}

\newcommand{\sym}{{\rm sym}}

\newcommand{\TV}{{\rm TV}}

\newcommand{\fnq}{{\rm fnq}}

\newcommand{\CF}{{\rm CF}}

\newcommand{\JEP}{{\rm JEP}}

\newcommand{\IND}{{\rm IND}}

\newcommand{\aut}{{\rm aut}}

\newcommand{\bd}{{\rm bd}}

\newcommand{\id}{{\rm id}}

\newcommand{\LST}{{\rm LST}}

\newcommand{\grp}{{\rm grp}}

\newcommand{\exlf}{{\rm exlf}}
\newcommand{\slf}{{\rm spc}}

\newcommand{\Rang}{{\rm Rang}}

\def\emptycs{}
\def\evaluateROMANlist{%
        \ifx\ROMANlist\emptycs\else
        \expandafter\xxevaluate\ROMANlist\xxfertig\evaluateROMANlist\fi}
\def\xxevaluate#1,#2\xxfertig{\expandafter\newcommand\csname#1\endcsname{\mathrm
{#1}}\def\ROMANlist{#2}}
 \newcommand{\ROMANlist}{Prop,pseudoProd,IDC,orth,xyz,}
 \evaluateROMANlist

\newcommand{\rest}{{\restriction}}

\newcommand{\then}{{\underline{then}}}
\newcommand{\when}{{\underline{when}}}

\newcommand{\If}{{\underline{if}}}
\newcommand{\Iff}{{\underline{iff}}}

\newcommand{\lqq}{{{`} {`} }}

\newcommand{\bfc}{{\mathbf c}}

\newcommand{\bfT}{{\mathbf T}}
\newcommand{\bfK}{{\mathbf K}}

\newcommand{\cE}{{\mathscr E}}

\newcommand{\cF}{{\mathscr F}}
\newcommand{\cG}{{\mathscr G}}

\newcommand{\gk}{{\mathfrak k}}

\newcommand{\cf}{{\rm cf}}

\newcount\skewfactor
\def\mathunderaccent#1#2 {\let\theaccent#1\skewfactor#2
\mathpalette\putaccentunder}
\def\putaccentunder#1#2{\oalign{$#1#2$\crcr\hidewidth
\vbox to.2ex{\hbox{$#1\skew\skewfactor\theaccent{}$}\vss}\hidewidth}}

\newbox\noforkbox \newdimen\forklinewidth
\setbox0\hbox{$\textstyle\bigcup$}
\setbox1\hbox to \wd0{\hfil\vrule width \forklinewidth depth \dp0
                        height \ht0 \hfil}
\setbox\noforkbox\hbox{\box1\box0\relax}
\def\unionstick{\mathop{\copy\noforkbox}\limits}
\def\nonfork#1#2_#3{#1\unionstick_{\textstyle #3}#2}
\def\nonforkin#1#2_#3^#4{#1\unionstick_{\textstyle #3}^{\textstyle
    #4}#2}
\setbox0\hbox{$\textstyle\bigcup$}
\setbox1\hbox to \wd0{\hfil{\sl /\/}\hfil}
\setbox2\hbox to \wd0{\hfil\vrule height \ht0 depth \dp0 width
                                \forklinewidth\hfil}
\newbox\doesforkbox
\setbox\doesforkbox\hbox{\box1\box0\relax}
\def\nunionstick{\mathop{\copy\doesforkbox}\limits}

\def\fork#1#2_#3{#1\nunionstick_{\textstyle #3}#2}
\def\forkin#1#2_#3^#4{#1\nunionstick_{\textstyle #3}^{\textstyle
    #4}#2}

\makeatletter
\@namedef{subjclassname@2020}{%
  \textup{2020} Mathematics Subject Classification}
\makeatother

\newcommand{\stickT}{%
\setbox255=\hbox{\raise1ex\hbox{$\hspace{0.2pt}\,\bullet\,$}}
\mathord{\rlap{\hbox to\wd255{\hss\hbox{$|$}\hss}}
\box255}
}
\newcommand{\stickS}{%
\setbox255=\hbox{\raise0.6ex\hbox{$\scriptstyle\bullet$}}
\mathord{\rlap{\hbox to\wd255{\hss\hbox{$\scriptstyle|$}\hss}}
\box255}
}

\newenvironment{PROOF}[2][\proofname.]
   {\begin{proof}[#1]\renewcommand{\qedsymbol}{\textsquare$_{\rm #2}$}}
   {\end{proof}}

\usepackage[hidelinks]{hyperref}
\begin{document}
\makeatletter\def\shfiuwefootnote{\gdef\@thefnmark{}\@footnotetext}\makeatother\shfiuwefootnote{Version 2023-02-24. See \url{https://shelah.logic.at/papers/1175/} for possible updates.}

\title{Canonical universal locally finite groups \\
Sh1175}
\author {Saharon Shelah}
\address{Einstein Institute of Mathematics\\
Edmond J. Safra Campus, Givat Ram\\
The Hebrew University of Jerusalem\\
Jerusalem, 9190401, Israel\\
 and \\
 Department of Mathematics\\
 Hill Center - Busch Campus \\ 
 Rutgers, The State University of New Jersey \\
 110 Frelinghuysen Road \\
 Piscataway, NJ 08854-8019 USA}
 
\email{shelah@math.huji.ac.il}

\urladdr{http://shelah.logic.at}

\thanks{The author thanks Alice Leonhardt for the beautiful typing of earlier versions (up to 2019) and, in later versions, the  author would like to thank the typist for his work and is also grateful for generous funding of typing services donated by a person who wishes to remain anonymous. The author thanks the Israel Science Foundation (ISF) (2019-2023) for partially supporting this research by grant 1838/19.  First typed February 18, 2016, as part of \cite{Sh:1098}.  In References \cite[0.22=Lz19]{Sh:312} means \cite[0.22]{Sh:312} has label z19 there, L stands for label; so will help if \cite{Sh:312} will change.   Also, \cite[Th.3.5=Th.1.5=Lb24]{Sh:1181} refer to Th.3.5 in the published version, Th. 1.5 in the arXive version, and label b24 in the latex file. The reader should note that the version in my website is usually more updated than the one in the mathematical archive}








\date{January 25, 2023}

\begin{abstract}
    We prove that for $\lambda = \beth_\omega$ or just $\lambda$
    strong limit singular of cofinality $\aleph_0$, if there is a
    universal member in the class $\mathbf K^{\rm lf}_\lambda$
    of locally finite groups of cardinality $\lambda$, then there is a canonical one
    (parallel to special models for elementary classes, which is the
    replacement of universal homogeneous ones and saturated ones in
    cardinals $\lambda = \lambda^{< \lambda}$).

    For this, we rely on the existence of enough  indecomposable such groups,
    as proved in $\lqq.$ Density of indecomposable locally finite groups". 
    We also more generally deal with the existence of universal members 
    in general classes 
    for such cardinals. 
\end{abstract}

\maketitle
\numberwithin{equation}{section}
\setcounter{section}{-1}
\newpage


%

%
%

%
%

%
%


\section{Introduction} \label{0}

\subsection{Background and aims}\ 

Our motivation is investigating 
the class 
$ \mathbf{K} _{\lf}$  
of locally finite groups (that is, groups such that every finitely generated sub-group is finite), so the reader
may consider only this case ignoring the general case; or consider universal classes (see Def. \ref{a29}).  
We continue \cite{Sh:312}, see history there and see more in \cite{Sh:1098}, and on earlier history 
the book \cite{KeWe73}.

We wonder:
\begin{problem}\label{w6}
1) Is there a universal $G \in \mathbf K^{\lf}_\lambda$
(= the class of members of $ \mathbf{K} _{\lf}$  
of cardinality $ \lambda $), 
see \ref{a28}(1);
e.g. for $\lambda =
\beth_\omega$?  Or just $\lambda$ strong limit of cofinality $\aleph_0$
(which is not above a compact cardinal)?

2) May there be a universal $G \in \mathbf K^{\lf}_\lambda$, when $\lambda < \lambda^{\aleph_0}$, e.g. for $\lambda =
\aleph_1 < 2^{\aleph_0}$, i.e. consistently?
\end{problem}

Generally, on the problem of the existence of a universal  model of a class in cardinality $ \lambda $ see  the classical Jonsson 
\cite{Jo56}, \cite{Jo60}, 
Morley-Vaught  \cite{MoVa62} 
and the recent survey \cite{Dj05} and later \cite{Sh:1151}.

Returning to locally finite groups, 
concerning \ref{w6}(1) recall that by Grossberg-Shelah \cite{Sh:174},  if $ \lambda = \lambda ^ {\aleph_0} $ 
then there is no universal member  
of $ \mathbf{K} ^{\lf}_\lambda $.  
\underline{But}   
if $\lambda$, 
a  
strong limit 
cardinal 
of cofinality, $\aleph_0$ is above a compact cardinal
$\kappa$, \then \, there is 
$G \in \mathbf{K} ^{\lf}_\lambda$   
which is universal.  
 So Problem \ref{w6} address the main open cases.

 More fully by \cite[2.17 = L s56]{Sh:1029} we have: 

\begin{conclusion}\label{s56}
    1) $\mu = \cf(\mu), \, \mu^{+} < \lambda = \cf(\lambda) < 2^{\mu}$ \underline{then} there is no group of cardinality $\lambda$ universal for the class of locally finite groups. 

    2) For example, if $\aleph_{2} \leq \lambda = \cf(\lambda) < 2^{\aleph_{0}}$ this applies.
\end{conclusion}

So the only other cases left are:

\begin{enumerate}
    \item[(A)] $\cf(\lambda) < \lambda < 2^{\aleph_{0}},$

    \item[(B)] $\beth_{\delta} < \lambda < 2^{\beth_{\delta}}, \, \cf(\delta) = \aleph_{0}$ and $\lambda$ singular. 
\end{enumerate} 

\underline{Question}: Are there partial long orders on $\{ {}^{\theta} G \colon G \in \bfK_{\lf} \}$?

Above, if $\aleph_{0} \, {<} \, \mu = \cf(\mu) < \lambda, \, 2^{\mu} > \lambda$ such that $\bfT_{J_{\mu}^{\bd}}(\lambda) < 2^{\mu},$ then there is no universal in $\lambda,$  as in \cite{Sh:409}, see \cite{Sh:F2150}.

Let us consider the model theory of locally finite groups.

Recall

\begin{definition}\label{x2}  
    1) $G$ is a lf \underline{(locally finite) group} \If \, $G$ is a group and  every finitely generated subgroup is finite. 
    
    2) $G$ is an exlf (\underline{existentially closed} lf) group  (in \cite{KeWe73} it
    is called ulf, universal locally finite group) \when \, $G$ is a locally finite group and for any finite groups  $K \subseteq L$ and embedding of $K$ into $G$, the embedding can be extended to an embedding of $L$ into $G$.
    
    3) Let $\mathbf K_{\lf}$ be the class of lf (locally finite) groups (partially ordered by $\subseteq$, being a subgroup) and let
    $\mathbf K_{\exlf}$ be the class of existentially closed $G \in \mathbf K_{\lf}$.
\end{definition}

Wehrfritz asked about the categoricity of the class of exlf groups  in any $\lambda > \aleph_0$.  This was answered by 
Macintyre-Shelah \cite{Sh:55} which proved that in every $\lambda > \aleph_0$ there are $2^\lambda$ non-isomorphic members of $\mathbf K^{\text{exlf}}_\lambda$.  This was disappointing in some sense: in $\aleph_0$ the class is categorical, so the question was perhaps
motivated by the hope that also general structures in the class can be understood to some extent.

The existence of a universal can be considered as a weak positive answer.  

A natural and frequent question on a class of structures is the existence of rigid members, i.e. ones with no non-trivial automorphism.  Now any exlf group $G \in \mathbf K_{\exlf}$ has non-trivial automorphisms - the inner automorphisms (recalling it has a trivial center).
So the natural question is about complete members where a group is called complete \Iff \, it has no non-inner automorphism.

Concerning the existence of a complete, locally finite group of cardinality $\lambda$: Hickin \cite{Hi78} proved one exists  in $\aleph_1$ (and more, e.g. he finds a family of $2^{\aleph_1}$ 
such groups pairwise far apart, i.e. no
uncountable group is embeddable into two of them).  Thomas \cite{Th86} assumed G.C.H. and built one in every successor cardinal (and
more, e.g. it has no Abelian or just solvable subgroup of the same cardinality).  Related are Giorgetta-Shelah \cite{Sh:83}, Shelah-Zigler \cite{Sh:96}, which investigate $\mathbf K_{G_\ast}$ getting similar results.

Dugas-G\"obel \cite[Th.2]{DgGb93} prove that for $\lambda = \lambda^{\aleph_0}$ and $G_0 \in \mathbf K^{\lf}_{\le \lambda}$ there is a
complete $G \in \mathbf K^{\exlf}_{\lambda^+}$ extending $G_0$; moreover $2^{\lambda^+}$ pairwise non-isomorphic ones.  Then Braun-G\"obel \cite{BrGb03} got better results for complete locally finite $p$-groups.

Now \cite{Sh:312} show that though the class $\mathbf K_{\text{exlf}}$ is very ``unstable"
there is a large enough set of definable types so we can imitate stability theory and have reasonable control in building exlf groups, using
quantifier-free types.  This may be considered a ``correction" to the non-structure results discussed above. This was applied to build 
a canonical extension of a locally finite group of the same cardinality which is existentially closed (it was known to exist in the power set, see \cite{KeWe73}). Also, there are  endo-rigid locally finite groups in a more relaxed way.


Returning to the present work,   here we deal with the universality problem for  $\mu = \beth_\omega$ or just strong limit of
cofinality $\aleph_0$.  We prove for $\mathbf K_{\lf}$ and similar classes that if there is a universal model of cardinality $\mu$, \then
\, there is something like a special model of cardinality $\mu$, in particular, universal, and unique up 
to isomorphism.  This relies on \cite{Sh:1181}, which proves the  existence and even   the  density of so-called
$\theta$-indecomposable (i.e. $\theta$ is not a possible cofinality) models in $\mathbf K_{\lf}$ of various cardinalities continuing
Carson-Shelah \cite{Sh:1169} 
which deal with the class of groups.

Returning to Question \ref{w6}(1), a possible avenue is to try to prove  the existence of universal members in $\mu$ when $\mu =
\Sigma_{n < \omega } \mu_n$ each $\mu_n$ measurable  $< \mu$, i.e. maybe for some
reasonable classes this holds.

We thank the referee for helpful remarks and later Mark Po\'or.  

\subsection{Definitions}\

\begin{context}\label{a27}

    $ \mathbf{K}$ will be one of the following cases:  
    
    \underline{Case 1}: $ \mathbf{K} = \mathbf{K} _{\lf}$, the class of locally 
    finite groups, so the submodel relation is just a subgroup,
    
    \underline{Case 2}: $ \mathbf{K} $ is a universal class, see 
    Def \ref{a29}(1)   
    below, the submodel 
    relation means just a submodel,
    
    \underline{Case 3} $ \mathbf{K} $ is 
     $\gk = (K_{\gk},\le_{\gk})$  an a.e.c.
     with  $ \LST_ \mathfrak{k} < \mu $ , see \cite[\S1]{Sh:88r};  we shall
     only comment on it. In particular, in this context,  in the definitions, 
      $ M \subseteq N $ should be replaced by 
      $ M \le _ \mathfrak{k} N $.
\end{context}

\begin{definition}\label{a28}
    1)  We say that $ M \in \mathbf{K}_ \mu $ is universal  (in $ \mathbf{K} $ or in $ \mathbf{K} _ \mu$, see \ref{a29}) \when \, every member of $ \mathbf{K} _{\mu }$ 
    can be embedded into it.
    
    2) We say that $ M \in \mathbf{K} $  is 
    universal for $ \mathbf{K} _{< \mu }$   
    \when \, 
    every $ M \in \mathbf{K}_{< \mu }$ can be 
    embedded into it; see Def \ref{a29}(4)  
    below.
    
    3) We define similarly   $\lqq $$M \in \mathbf{K} $ is universal for $\mathbf{K} _{\mu }$"  and $\lqq $$ M \in \mathbf{K} $ is universal for $ \mathbf{K} _{\le \mu }$".
    \end{definition} 
    
\begin{definition}\label{a29}
    1) We shall say that $ \mathbf{K} $  is a universal class
    \when \,
    for some vocabulary $ \tau = \tau _ \mathbf{K}$:
    
    \begin{enumerate} 
    \item[(a)]
    $ \mathbf{K} $  is a class of $ \tau $-models, 
        closed under isomorphisms,
    \item[(b)] a $ \tau $-model belongs to $ \mathbf{K} $ 
     iff every 
    finitely generated sub-model belongs to it,
    \end{enumerate}

    3) Let $ \mathbf{K} _\mu  $ be the class of $ M \in \mathbf{K} $ 
    of cardinality $ \mu  $. We define 
    $ \mathbf{K} _{< \mu }, \mathbf{K} _{\le \mu }$ 
    naturally.  

    4) For cardinals $ \lambda \le  \mu $ let $ \mathbf{K} _{\mu, \lambda }$
    be the class of pairs $ (N,M)$  such that $ N \in \mathbf{K} _ \mu  ,
    M \in \mathbf{K} _ \lambda $ and $ M \subseteq N$.
    
    5) Let $ ( N_1, M_1) \le _{ \mu, \lambda  }(N_2, M_2 ) $ mean that   
    $ (N_ {\ell} , M_ {\ell} ) \in \mathbf{K} _{ \mu, \lambda  }$ 

    for $ {\ell} = 1,2 $ and
    $ M_1 \subseteq M_2, N_1 \subseteq N_2$.

     6) For $\lambda \leq \mu$ we define $ \mathbf{K} _{\mu, < \lambda }$
     and $ \le _{\mu, < \lambda }$ similarly.
     
     7) A universal class $ \mathbf{K} $ can be considered  as the  a.e.c. $ \mathfrak{k} = (\mathbf{K}, \subseteq)$
\end{definition} 

\begin{notation}\label{a30}
    1) Let $M,N $ and also $ G,H,L $ denote members of $ \mathbf{K}$.
    
    2) Let $ |M| $ be the universe = set of elements of $ M $  and $ \| M \| $ its cardinality.
    
    3) Let $ a, b, c ,  d $ denote members of such $ M$, and $ \bar{ a }, \bar{ b } \dots$  denote sequences of such elements.
\end{notation}

\begin{definition}\label{a32}
    1) We say the pair $(N, M)$ is an $(\chi,\mu,\kappa)$-amalgamation
    base (or amalgamation pair, but we may omit $ \chi $ when
    $ \chi = \mu$, and we may even omit $ \mu, \kappa $  too when clear from the context) \when \,:  
    \begin{enumerate} 
        \item[(a)] $ (N,M ) \in \mathbf{K} _{\mu, \kappa },$
        
        \item[(b)] if $ N_1=N $ and $ M \subseteq  N_2 
            \in \mathbf K_{  
            \chi} $
        
        \then \, 
        $ N_1, N_2 $ can be amalgamated over $ M$, this mean that 
        for some
        $N_3,f_1,f_2$ we have $M \subseteq  N_3 \in \mathbf K$ and
        $f_\ell $-embeds $N_\ell$ into $N_3$ over $M$ for $\ell = 1, 2$.
    \end{enumerate}

     2) We say that the   pair   $ (N,M)$ is a universal   $ (\mu, \lambda )$-amalgamation 
    base 
    (we may omit $ \mu, \lambda $)   
    \when \,:
     \begin{enumerate} 
     \item[(a)] $ (N,M) \in \mathbf{K} _ {\mu , \lambda },$
     
     \item[(b)] if $ N \subseteq N' \in \mathbf{K} _\mu $ 
     \then \, $ N' $ can be embedded into $ N $ over $ M $.
    \end{enumerate} 
     
    3) We  may in  parts (1),(2) omit  $  \mu , \kappa $ when $ ( \mu, \lambda )=(\| N \|, \| M \|)$. 
    
    4) We say $M \in \bfK_{< \mu}$ is an amalgamation base inside $\bfK_{< \mu}$ when: if $M \leq N_{\ell} \in \bfK_{< \mu},$ then $N_{1}, N_{2}$ can be amalgamated over $M$ (see \ref{a32}(1)(b)) but $N_{3} \in \bfK_{< \mu}.$ 
\end{definition}

\newpage

\section{Indecomposability}\label{1}

In this section we deal with indecomposability, equivalently $\CF(M)$, see e.g. \cite{Sh:524}; we have $\bfK_{\lf}$ in mind but still it is meaningful and of interest also for other classes.

Why do we deal with indecomposable 
members $ \mathbf{K} $? 
When we shall try to understand universal 
members $ M $ of $ \mathbf{K} _\mu $ we  
shall use some $ \theta $-indecomposable 
$ N \subseteq M $ of cardinality $ < \mu $. 
How will  this   help us? The point is that 
$ N \in \mathbf{K} _{< \mu }$ may have too many 
embeddings into $ M $, but if 
$\theta = \cf( \theta ) \neq \cf(\mu ),$    
$ \alpha < \mu \Rightarrow | \alpha | ^{\|  N\|} < \mu,$  
$ N $  is $ \theta $-indecomposable
and $ \theta $ is regular uncountable $ < \mu,$ then this is not the case. 

We need indecomposable $ \mathbf{c} :[ \lambda] \rightarrow \theta $  
in order to build enough $ \theta $-indecomposable locally finite groups
(as done in \cite{Sh:1181}).  


\begin{definition}\label{b8}
    1) We say $M$ is $\theta$-indecomposable  or $ \theta \in \CF(M)$  
    \when \,: $\theta$ is regular
    and if $\langle M_i:i < \theta\rangle$ is $\subseteq$-increasing with
    union $M$, \then \, $M = M_i$ for some $i$.
    
    2) We say $M$ is $\Theta$-indecomposable \when \, it is
    $\theta$-indecomposable for every $\theta \in \Theta$.
     We say $M$ is $\Theta^{\orth}$-indecomposable \when \, it is
    $\theta$-indecomposable for every  regular  $\theta \notin \Theta$.

    3) We say $G$ is $\theta$-indecomposable inside $G^+$ \when \,:
    
    \begin{enumerate}
    \item[(a)]  $\theta = \cf(\theta)$;
    
    \item[(b)]  $G \subseteq G^+$;
    
    \item[(c)]  if $\langle G_i:i \le \theta\rangle$ is
      $\subseteq$-increasing continuous and  
      $ G_\theta =G^+$ 
        (hence $ G \subseteq G_ \theta $)  
     \underline{then} \, for some $i < \theta$ we have $G \subseteq G_i$.
    \end{enumerate}

    4) For $ \theta = \cf(\theta )  \le  \lambda \le \mu $  such that $ \theta \notin \Theta _ \lambda$ (see \ref{b24}(1)) we say $ \mathbf{K} $ is $ (\mu,  \lambda, \theta )$-indecomposable \when \,  for every pair $ (N,M )\in \mathbf{K} _{\mu, \lambda }$ there is $ (N_1, M_1 ) \in \mathbf{K} _{\mu, \lambda }$
    which is $  \le _{\mu, \lambda }$-above it and $ M_1 $  is $\theta $-indecomposable  (really, not just  
    inside $ N_1 $).  For $ \theta = \cf(\theta )  <   \lambda \le \mu $
    we   say  $ \mathbf{K} $ is  
    $ (\mu, < \lambda , \theta ) $-indecomposable  \when: 
    
    if $ \theta = \cf( \theta ) \le \lambda _1 < \lambda , \theta \notin   
        \Theta _{\lambda _1}$ then $ \mathbf{K} $  is 
        $ (\mu , \lambda _2, \theta )$-indecomposable for some $\lambda_{2} \in [\lambda_{1}, \mu).$ 

    5) We say $\mathbf c:[\lambda]^2 \rightarrow S$ is
    $\theta$-indecomposable \when \,: if $\langle u_i:i < \theta\rangle$ is
    $\subseteq$-increasing 
        sequence of sets 
    with union $\lambda$ then $S = 
    \{\mathbf c\{\alpha,\beta\}:\alpha \ne \beta \in u_i\}$ for some $i <
    \theta$; 
    
    6) We may replace above the cardinal $ \theta $  by a set or class $ \Theta $ of regular cardinals, 
    (as done 
    in \ref{b8}(2)).  
\end{definition}

 A group $G$ may be 
 considered 
 indecomposable as a group or as a semi-group;
our default choice is semi-group; but note that for locally finite
groups the two 
interpretations are equivalent.

The following was proved in \cite{Sh:1181}.  

\begin{theorem} \label{b24}
1)
 If $\lambda \ge \aleph_1$ and we let 
$ \Theta_\lambda  = \{ \cf(\lambda )\} $ except that 
$ \Theta_\lambda  = \{ \cf( \lambda ), \partial \} 
    = \{\lambda, \partial  \} $
when $ (c)_{\lambda, \partial }$ below holds, 
\then \,  
clauses 
(a),(b) hold, where:

\begin{enumerate}
\item[(a)]  some $\bfc \colon [\lambda]^2 \rightarrow \lambda$ is
  $\theta$-indecomposable  for every  $ \theta = 
  \cf  (  \theta ) \notin \Theta_\lambda  $

\item[(b)]  for every $G_1 \in \bfK^{\lf}_{\le \lambda}$ there is an
extension $G_2 \in \bfK^{\lf}_\lambda$ which is 
$ \Theta ^{\orth}_\lambda $-indecomposable

\item[(c)$_{\lambda, \partial  }$]  for some
$\mu,\lambda =\mu^+,\mu > \partial  = \cf(\mu)$ and
  $\mu = \sup\{\theta  < \mu:\theta $
  is a regular Jonsson cardinal$\}$.
\end{enumerate}

2) If $ \mu \ge \lambda \ge   
    \theta = \cf(\theta )$ ,   
    and $ \theta \notin \Theta _{\lambda }, \lambda \geq {\aleph_1} $  
\then \, 
$ \mathbf{K}_{\lf } $ is $ (\mu ,\lambda, \theta )$-indecomposable.

2A) In fact (on part (2)) it suffice to assume $(\exists \lambda_{1})(\lambda \leq \lambda_{1} \leq \mu \wedge \theta \notin \Theta_{\lambda_{1}}).$

3)   
If $ \mu \ge \lambda $  and $ (H_1, G_1 )\in \bfK_{\le \mu, \le \lambda }$
    \then  \,  
we can find a pair $ (H_2, G_2 )  \in K_{\mu ,\lambda }$  such that:
\begin{enumerate}
    \item[(a)]   $ G_2 $ is $ \Theta ^{\orth}_\lambda $-indecomposable,
    
    \item[(b)]  if $ \mu > \lambda$ then  $G_1$ is $\theta $-indecomposable inside $H_{2}$ for every regular $ \theta,$ 
    
    \item[(c)]  $H_2$ is $\Theta^\orth _\mu  $-indecomposable.  
\end{enumerate} 
\end{theorem} 

\begin{PROOF}{\ref{b24}}
1) By\footnote{But Theorem 1.5 in the author's archive version, similarly 3.4  is 1.4} \cite[Th.~3.5]{Sh:1181}.

2), 2A) The proof will serve also for part (3). 
Let $ (N,M) \in \mathbf{K} _{\mu, \lambda }$
be given.  
We choose a pair $ (\chi , \partial )$ of cardinals 
    and $ \mathbf{c} $   
such that 
$ \lambda \le \chi \le \mu $
  $ \partial = \cf( \partial ) \le \lambda,
        \partial \not= \theta $ and 
$ \mathbf{c} \colon [\chi ]^2 \rightarrow   \chi $  
        is $ \theta $-indecomposable; 
        (possible here as $ \theta \notin \Theta _ \lambda , 
            \lambda \geq {\aleph_1} $ even for $ \chi = \lambda$).  

 By induction of $ \alpha \le \partial   $ we choose
 $ H_ \alpha, L_\alpha $, 
but $ L_ \alpha $ is chosen together with $ H_{\alpha + 1 }$ 
    when $ \alpha $ is a successor ordinal,  
 such that:
 
 \begin{enumerate} 
     \item[(a)] $ (H_\alpha, L_ \alpha ) \in \mathbf{K} _{\mu, \lambda   }$ 
     is increasing continuous with  $ \alpha $ 
     \item[(b)] $ (H_0, L_0 )= (N,M),$
     
     \item[(c)] if  $\alpha =  \beta  + 1 < \theta $   then 
     and $ L_{\alpha}  $ is 
     $ \theta $-indecomposable.
 \end{enumerate} 
 
Why can we carry the induction?
For $ \alpha = 0 $ this is trivial; 
similarly for $ \alpha $
a limit ordinal.
Lastly by clause (b) of part (1), for $ \alpha = \beta + 1 \le \alpha _\ast $,
recall the proof of \cite[3.4]{Sh:1181}, pedantically  as without loss of generality $H_{\beta}, \, L_{\beta}$ are existentially closed hence generated by the elements of order $2,$ let $\langle a_{\alpha} \colon \alpha < \mu \rangle$  list $\{ a \in L_{\beta} \colon a \text{ of order } 2 \}.$ By \cite[Prop. 3.4(2)]{Sh:1181}, with $u_{\alpha} = \{ \alpha \},$ we can find $H_{\alpha, 1} \in \bfK_{\mu}^{\mathrm{lf}}$ extending $H_{B}$ and pairwise commuting $b_{\alpha} \in H_{\alpha, 1}$ each of order $2,$ for $\alpha < \mu$ (the order $2$ was not mention but proved) and pairwise commuting $d_{\alpha} \in H_{\alpha, 1},$ each of order $2,$ for $\alpha < \mu$ such that, $L_{\beta}$ is included in the subgroup $L_{\alpha, 1}$ of $H_{\alpha, 1}$ generated by $\{ b_{\alpha}, d_{\alpha} \colon \alpha < \lambda \}.$ 

Now apply \cite[Prop. 3.4(1)]{Sh:1181}
 for a $\theta$-indecomposable $\bfc \colon [\lambda]^{2} \to \lambda.$


3) We deal with every regular $ \theta  \le \mu  $
successively. Fixing $ \theta $ we can 
    use 
the proof 
of part (2). 
\end{PROOF} 


Now comes the central definition, 
what is its role? 

We like to sort out when there is a 
universal member of $ \mathbf{K} _ \mu $ 
and when there is a canonical universal 
member. For reasons explained above we 
concentrate on 
the case $ \mu $ is strong limit
of cofinality $ {\aleph_0} $, for example 
$ \beth _ \omega $. 
    To find out  the answer to those two questions  for every  universal class $ \mathbf{K} $ seem too much to hope for. The Def \ref{b28} accomplishes a more modest task: it gives a large frame satisfied  by a large family of pairs $ (\mathbf{K}, \mu )$   for which we shall prove an equivalence.  In particular our class $ \mathbf{K} _{\lf}$ belongs  to this family.


\begin{definition} \label{b28} 
We say that $ \mathbf{K} $ is $ \mu $-nice \when \,: 

\begin{enumerate}
    \item[(a)] $ \tau _ \mathbf{K} $ has cardinality $ < \mu $, 
    
    \item[(b)]  for every $ M \in \mathbf{K} _{< \mu }$ there is
    $ N \in \mathbf{K} _ \mu $ extending $ M $, 
    
    \item[(c)] $ \mathbf{K}  $ has the JEP (joint embedding property),
    
    \item[(d)] $ \mathbf{K} $ is  $ ( \mu, < \mu   
    , \cf( \mu ))$-indecomposable or just,

    \item[(d)$'$] for arbitrarily large $ \lambda _2 < \mu $  letting $ \theta = \cf( \mu )  \le \lambda _2 $  
    we have $ \mathbf{K} $ is $(\mu, \lambda _2, \theta  )$-indecomposable. 
\end{enumerate}
\end{definition} 

Naturally we like to prove that   the pair   $ (\mathbf{K} _{\lf}, \beth _ \omega )$ falls under the frame of Def \ref{b28}.  This is the role of \ref{b31}, \ref{b34}.  In \S3   we point out an additional family. For the main case, $ \mu $ is a strong limit of cofinality $ {\aleph_0} $. 

\begin{claim}\label{b31}
     $ \mathbf{K} _{\lf}$ is $ \mu$  -nice when $ \mu \geq {\aleph_1}.$ 
\end{claim} 

\begin{PROOF}{\ref{b31}}
In Def \ref{b28} clause (a)  
is trivial and as $ \mathbf{K} _{\lf}$  is closed 
under products clearly  clauses 
(b),(c) are clear and clause (d),  
    for $ \mu $  
    regular is trivial (and is not used)
    and for $ \mu $   
    singular it   
holds by \ref{b24}(3), see also    \ref{b34}(2)  below.
\end{PROOF} 


We give below more than what is strictly needed.
\begin{claim} \label{b34}

Assume $ \mathbf{K}= \mathbf{K} _{\lf} $.

1)
We have $ (A) \Rightarrow (B) $ where:
\begin{enumerate} 
\item[(A)] 
\begin{enumerate} 
\item  $  \mu \ge {\aleph_1},$  

\item $ \delta  _\ast \le \mu  $ and $ \lambda_ \alpha 
< \mu$ 
for $ \alpha < \delta _\ast,$

\item  $ \lambda _\alpha  \ge | \alpha |  $  is non-decreasing, 



\item $ G_ 1 \in \mathbf{K} _{\le \mu },$

\item $ G_{1, \alpha} \in \mathbf{K} _{\le \lambda _\alpha}$ is a subgroup of $G_{1}$
for $ \alpha < \delta _\ast,$ 
\end{enumerate}

\item[(B)] There are $ G_2, \bar{G}_2 $ such that:
\begin{enumerate} 
\item $ G_2 \in \mathbf{K} _{\mu} $  extends $ G_1,$

\item $ \bar{G}_2 = \langle G_{2, \alpha } \colon \alpha <  \delta  _\ast\rangle    $ 
    is increasing, 
    
\item $ G_{2, \alpha } \in \mathbf{K} _{\lambda _\alpha }$ 
 extend $G_{1,\alpha },$
 
\item $ G_ 2 $ is $ \Theta $-indecomposable
where\footnote{If $\mu = \sum \{ \lambda_{\alpha} \colon \alpha < \delta_{\ast} \}$ then $\cf(\mu) = \cf(\delta_{\ast}),$ hence the ``$\cup \{ \delta_{\ast} \}$'' is not necessary.} 
$ \Theta =  (\Theta_{\mu} \cup \{ \cf( \delta _\ast)\})  ^{\orth}$
\item $ G_{2, \alpha } $ is 
$ \Theta ^{\orth}_{\lambda_\alpha } $-indecomposable
(not just inside $ H_2)$
for every $ \alpha < \delta  _\ast $

\item  if $ \mu =\Sigma \{\lambda _ \alpha \colon \alpha < \delta _\ast \}$ 
 the $ G_2 = \cup  \{ G_{2, \alpha } \colon \alpha < \delta _\ast  \},$

 \item if $\mu > \sum \{ \lambda_{\alpha} \colon \alpha < \delta_{\ast} \}$ then $G_{2}$ is $\Theta_{M}^{\rm{orth}}$-indecomposable.
\end{enumerate} 
\end{enumerate}

 2) If $ \mu > \lambda \ge {\aleph_1} $ \then \,  $ {\aleph_0}  \in \Theta  ^{\orth}_{\cf(\mu)}  \cup   \Theta ^{\orth}_\lambda $  except possibly when   
   $ \mu = \lambda  ^+, \cf( \lambda  )= {\aleph_0} $.   
\end{claim} 


\begin{PROOF} {\ref{b34}}
1)
 By induction of $ \alpha \le \delta _\ast $ we choose
 $ H_ \alpha,\bar{ H}_ \alpha , 
    L_\alpha $, but $ L_ \alpha $ is chosen together
 with $ H_{\alpha + 1 }$ and not chosen for $ \alpha = \alpha _\ast $, 
 such that:    
 
 \begin{enumerate} 
 \item[(a)] $ H_\alpha $ 
 is increasing continuous with 
 $ \alpha $ 
 \item[(b)] $ H_ 0 = G_1  $  
    and $ \alpha > 0 \Rightarrow H_ \alpha \in \mathbf{K}_{\mu} $ 
 \item[(c)] $ (H_\alpha , L_ \beta  )\in \mathbf{K} _{\lambda, \lambda _\beta  }$
   when $ \alpha = \beta + 1 \le  \alpha _\ast$
 \item[(d)]  $  \bar{ H}_ \alpha = 
    \langle H_{\alpha , \varepsilon }: \varepsilon < \delta _\ast\rangle $
   such that if $ \mu = \Sigma \{\lambda _ \varepsilon : \varepsilon < \delta _\ast \}$ then 
   this sequence is increasing with union $ H_ \alpha $ 
   and $ H_{\alpha, \varepsilon } $ has cardinality 
   $ \lambda _ \varepsilon $  
        when $ \alpha > 0 $ and $ \le \lambda _ \varepsilon $ when $ \alpha = 0 $
 \item[(e)] $ G _ {1, \beta },H_{ \beta , \varepsilon },  L_ \gamma$
 are sub-groups of $ L_ \alpha $ when $ \beta \le \alpha, 
 \varepsilon  \le \alpha , \gamma < \alpha $ 
 \item[(f)]  $ L_ \beta  $ is 
 $ \Theta ^{\orth}_{\lambda_\beta  }$-indecomposable,
 \item[(g)] $ G_2 $ is $ \Theta $-indecomposable  
   where $ \Theta = (\Theta _ \mu \cup \{ \cf (\delta _\ast)\} )^\orth$
 \end{enumerate} 
 
 Why can we carry the induction? 
 We choose $ \bar{ H } _ \alpha $ just after $ H_ \alpha $  was chosen. For $ \alpha = 0 $ this is trivial (note that  $ L_ \alpha $ is not chosen), similarly for $ \alpha $ a limit ordinal. Lastly  for  $ \alpha = \beta + 1 \le \alpha _\ast $, Definition  \ref{b8}(4) \ref{b24}(3)   gives the desired  conclusion.  In details, first choose $ L^+_ \beta \subseteq H_ \beta $ of cardinality at most $ \lambda _ \alpha $ containing  the desired sets (listed in clause (e)). Then apply \ref{b24}(3) to the pair $ (H_ \beta, L^+_ \beta )$  to get $ (H_\alpha, L_ \alpha  )$. Lastly, if $\mu > \sum \{ \lambda_{\alpha}^{+} \colon \alpha < \delta_{\ast} \}$, let $ G_2 \in \mathbf{K} _{\lambda }$ extend $ H_{\alpha_{\delta_{\ast}}}$ and satisfies the indecomposablity demand and, if $\mu > \sum \{ \lambda_{\alpha} \colon \alpha < \delta_{\ast} \},$ let $G_{2} = H_{\delta_{\ast}}.$ Now, letting
$ G_{2, \alpha} = L_ \alpha $ we are done.

2) Easy.
 \end{PROOF} 

\begin{claim}\label{b37}
    If $ \mu $ is strong limit singular and $ N \in \mathbf{K} _ \mu $ 
    \then \, the set $\IDC_{< \mu  }(N)$  
    has cardinality $ \le \mu$ where, for $ N \in \mathbf{K} _ \mu $, 
    
    \begin{enumerate} 
        \item[($\ast$)] 
        $ \IDC_{< \mu }(N)= \{M : M \subseteq N $ has cardinality $ < \mu $ 
        and is $ \cf(\mu )-$indecomposable $ \}  $.
    \end{enumerate} 
\end{claim}
 
\begin{PROOF} {\ref{b37}} 
    Easy. 
\end{PROOF} 

\newpage

\section{Universality}\label{2}

For quite many classes, there are universal members in any (large enough)
$\mu$ which is strong limit of cofinality $\aleph_0$, see \cite{Sh:820} which includes history.    Below we investigate ``is there a universal member of $\mathbf K^{\lf}_\mu$ for such $\mu$".  We prove that if there is a universal member, e.g. in $\mathbf K^{\lf}_\mu$, then there is a canonical one. 

What do we mean by $\lqq$canonical"?   This is not a precise definition, but we mean it is unique  up to isomorphism, by a natural definition.  Examples we have in mind are the algebraic closure of a  field,  
the saturated model of a complete first order theory $ T $ in cardinality 
$ \mu ^+ = 2^ \mu > |T|$ and the special model 
of a complete first order theory $ T $ in a singular  strong limit cardinal $ \mu > | T |$, see \cite{ChKe62}.  The last one means: 

\begin{enumerate} 
    \item[($\ast$)]  for such $ T, \mu $ we say that $ M $ is a special
    model of $ T $ of cardinality $ \mu $  when some  $ \bar{ M } $ witness it
    which means:
    \begin{enumerate} 
        \item[(a)] $ \bar{ M } = \langle M_i : i < \cf(\mu )\rangle,$
        
        \item[(b)]  $ M_ i $  is $\prec $-increasing with $ i,$
        
        \item[(c)] each $ M_ i $ has cardinality $ < \mu,$
        
        \item[(d)] $M = \cup \{M_i \colon i < \cf( \mu ) \},$
        
        \item[(e)] for every $ \lambda < \mu $ for every large enough
        $ i < \cf( \mu ) $ the model $ M_ i $ is $ \lambda ^+ $-saturated.
    \end{enumerate} 
\end{enumerate} 

Considering our main case, $ \mathbf{K} _{\lf }$,   a major difference  (between what we prove here (e.g. for $ \mathbf{K} _{\lf})$ and $ (\ast) $ is that  here amalgamation fail, so clause (B) of \ref{b40}  is a poor man  replacement.

\begin{theorem}\label{b40}
    Assume that $\mu$ is a strong limit of cofinality $\aleph_0$ and $ \mathbf{K} $ is $ \mu $-nice.
    
    1) The following conditions are equivalent:
    
    \begin{enumerate}
    \item[$(A)$]   There is a universal $G \in \mathbf K^{}_\mu.$
    
    \item[$(B)$]   If $H \in \mathbf K^{}_\lambda$ is  $\aleph_0$-indecomposable for some $\lambda < \mu$, \then \, there is a sequence $\bar G = \langle G_\alpha \colon \alpha < \alpha_\ast \le \mu\rangle$ such that:
    
    \begin{enumerate}
    \item[$(a)$]   $H \subseteq G_\alpha \in \mathbf K^{}_\mu,$
    
    \item[$(b)$]   if $G \in \mathbf K^{}_\mu$ extend $H$, \then \, for some $\alpha,G$ is embeddable into $G_\alpha$ over $H$.
    \end{enumerate}
    
    \item[$(B)^+$]  We can add in (B)
    
    \begin{enumerate}
    \item[$(c)$]   if $\alpha_1 < \alpha_2 < \alpha_\ast$, \then \,  $ G_{\alpha _1}, G_{\alpha _2}$ cannot be amalgamated over $ H $,  that is there are no $G,f_1,f_2$ such that $H \subseteq G \in \mathbf K_{}$ and $f_\ell$ embeds $G_{\alpha_\ell}$ into $G$ over $H$ for $\ell=1,2$,

    \item[$(d)$]   $(H,G_\alpha)$ is an amalgamation pair (see Definition \ref{a32}(1)),   moreover a universal amalgamation base (see \ref{a32}(2)).
    \end{enumerate}
    \end{enumerate}

    2)  We can add in part (1):
    
    \begin{enumerate}
    \item[$(C)$]   there is $G_\ast$ such that:
    
    \begin{enumerate}
    \item[$(a)$]  $G_\ast \in \mathbf K^{}_\mu$ is universal 
            for 
    $\mathbf K^{}_{<\mu}$;
    
    \item[$(b)$]  $\cE^{\aleph_0}_{G_\ast,<\mu}$, see  see  
        Def.  
        \ref{b42} below, 
    is an equivalence relation with $\le \mu$ equivalence classes;
    
    \item[$(c)$]  $G_\ast$ is $ \mu$-special, see  \ref{b42}(E) below.
    \end{enumerate} 
    \end{enumerate} 
    \begin{enumerate} 
    \item[$ (C)^+ $]  like clause (C) but we add
    \begin{enumerate} 
        \item[(d)]  If $ G , G_\ast\in \mathbf{K} _ \mu $ are  $ \mu $-special 
        \then \, $ G, G_\ast $ are isomorphic, (that is uniqueness).
    \end{enumerate}
    \end{enumerate}
\end{theorem}

Before we prove \ref{b40}, we define (this definition  is not just used in the proof but also in phrasing  \ref{b40}(2)).  

\begin{definition}\label{b42}
    For $\theta = \cf(\mu) < \mu$  and  $M_\ast \in \mathbf{K} _\mu$: we define:
    
    \begin{enumerate}
    \item[(A)]   $\IND^\theta_{M_\ast,<\mu} = \{N:N \le_{\gk} M_\ast$ has
      cardinality $< \mu$ and is $\theta$-indecomposable$\}$.
    
    \item[(B)]   $\cF^\theta_{M_\ast,<\mu} = \{f$: for some
      $\theta$-indecomposable $N = N_f \in K^{}_{< \mu}$ with universe
      an ordinal, $f$ is an embedding of $N$ into $M_\ast\}$.
    
    \item[(C)]   $\cE^\theta_{M_\ast,<\mu} = \{(f_1,f_2) \colon f_1,f_2 \in
      \cF^\theta_{M_\ast,<\mu},N_{f_1} = N_{f_2}$ and there are
      embeddings $g_1,g_2$ of $M_\ast$ into some
      extension $M  \in \mathbf{K} _ \mu $ 
        of $M_\ast$ such that $g_1 \circ f_1 = g_2
      \circ f_2\}$.
      
    \item[(D)]    We say 
    $M_\ast$ is $\theta-\cE^\theta_{M_\ast,<\mu}$-indecomposably homogeneous
    (or just $ M_\ast $ is $ \theta $-indecomposably homogeneous)
     \when \ some $\bar{M}$ witness it, which mean: 

     \begin{enumerate}
         \item[(a)] $\bar{M} = \langle M_{i} \colon i < \cf(\mu) \rangle$ is increasing continuous with limit $M,$ 
         
         \item[(b)] if $f_1,f_2 \in \cF^\theta_{M_\ast,<\mu}, \, (f_1,f_2) \in \cE^\theta_{M_\ast,<\mu}$ and $(\exists i < \theta)( A \subseteq M_{i}), A$ of cardinality $< \mu,$  \then \, there is $(g_1,g_2) \in \cE^\theta_{M_\ast,<\mu}$ such that $f_1 \subseteq g_1 \wedge f_2 \subseteq g_2$ and $A \subseteq \Rang(g_1) \cap \Rang(g_2)$; it follows that if $\cf(\mu) = \aleph_0$ then for some $g \in \aut(M_\ast)$ we have $f_2 = g \circ f_1$.
     \end{enumerate}
     
    \item[(E)]  We say
      that $ M_\ast \in \mathbf{K} _\mu $ is $ \mu $-special
      \when \, it is  $ \theta $-indecomposably homogeneous 
      and  
      is universal for $ \mathbf{K} _{< \mu }$, that is 
      every $ M \in \mathbf{K} _{< \mu}$ is embeddable into it.
    \end{enumerate}
\end{definition}

\begin{remark}\label{b44}
    We may consider in \ref{b40} also $(A)_0 \Rightarrow (A)$ where
    
    \begin{enumerate}
        \item[$(A)_0$]   if $\lambda < \mu,H \subseteq G_1 \in \mathbf
          K^{}_{< \mu}$ and $|H| \le \lambda$, \then \, for some $G_2$ we
          have $G_1 \subseteq G_2 \in \mathbf K^{}_{< \mu}$ and $(H,G_2)$ is a
          $(\mu,\mu,\lambda)$-amalgamation base.
    \end{enumerate}
\end{remark}

\begin{PROOF}   {\ref{b40}}
of \ref{b40}  

It suffices to prove the following implications:
\medskip

\underline{$(A) \Rightarrow (B)$}:

Let $G_\ast \in \mathbf K^{}_\mu$ be universal and choose a sequence
$\langle G^\ast_n:n < \omega\rangle$ such that $G_\ast =
\bigcup\limits_{n} G^\ast_n,G^\ast_n \subseteq G^\ast_{n+1},|G^\ast_n| < \mu$.

Let $H$ be as in \ref{b40}(B) and let
$\cG = \{g:g$ embed $H$ into $G^\ast_n$ for some $n\}$.  So clearly
$|\cG| \le \sum\limits_{n} |G^\ast_n|^{|H|} \le \sum\limits_{\lambda <
  \mu} 2^\lambda = \mu$, (an over-kill).

Let $\langle g^\ast_\alpha:\alpha < \alpha_\ast \le \mu\rangle$ list $\cG$
and let $(G_\alpha,g_\alpha)$ be such that (exist by renaming):

\begin{enumerate}
\item[$(\ast)_1$]  $(a) \quad H \subseteq G_\alpha \in \mathbf K^{}_\mu$;

\item[${{}}$]  $(b) \quad g_\alpha$ is an isomorphism from $G_\alpha$
  onto $G_\ast$ extending $g^\ast_\alpha$.
\end{enumerate}

It suffices to prove that $\bar G = \langle G_\alpha:\alpha <
\alpha_\ast\rangle$ is as required in 
    clause (B).   
Now clause (B)(a) holds by
$(\ast)_1(a)$ above.  As for clause (B)(b), let $G$ satisfy 
$H \subseteq G \in \mathbf
K^{}_{\le \mu}$, so there is $G' \in \bfK_{\mu}$ extending $G,$ hence we can find an embedding $g$ of $G'$ into
$G_\ast$.  We know that $g(H) \subseteq G = \bigcup\limits_{n} G_n$ hence $\langle
g(H) \cap G_n:n < \omega\rangle$ is $\subseteq$-increasing with union
$g(H)$; but $g(H)$ by the assumption on $H$ is
$\aleph_0$-indecomposable, hence $g(H) = g(H) \cap G^\ast_n 
\subseteq G^\ast_n$ for some $n$. 
This implies 
$g \rest H \in \cG$ 
and 
so for some $ \alpha < \alpha _\ast $   
we have $ g \, {\rest} \, H = g^\ast_\alpha $. Hence 
$ g_\alpha ^{-1}g $ is an embedding of $ G $ 
into $ G_\ast $ extending
$ (g_\alpha  \rest H)^{-1}(g \rest H ) =
(g^\ast_\alpha )^{-1}   
    (g^\ast_\alpha )=\id_H $  as promised.
\medskip

\underline{$(B) \Rightarrow (B)^+$}:

What about $(B)^+(c)$? while $\bar G$ does not necessarily satisfy it,
we can ``correct it", e.g. we choose $u_\alpha,v_\alpha$ and if
$\alpha \notin \cup\{v_\beta:\beta < \alpha\}$ also $G'_\alpha$ by
induction on $\alpha < \alpha_\ast$ such that (the idea is that if $\beta
\in v_\alpha$ 
then  $ \beta > \alpha $  and   
$ G_\beta$ is discarded being embeddable into some $G'_\alpha $ and
$G'_\alpha$ is the ``corrected" member):

\begin{enumerate}
\item[$(\ast)^2_\alpha$]
\begin{enumerate}
\item[(a)]  $G_\alpha \subseteq G'_\alpha \in \mathbf K^{}_\mu$ if
  $\alpha \notin \cup\{v_\beta:\beta < \alpha\}$;

\item[(b)]  $u_\alpha \subseteq \alpha$ and $v_\alpha
  \subseteq \alpha_\ast \backslash (\alpha  + 1)$;  

\item[(c)]  if $\beta < \alpha$ then $u_\beta =
  u_\alpha \cap \beta$ and $u_\alpha \cap v_\beta =\emptyset$;

\item[(d)]   if $\alpha = \beta +1$ then $\beta \in
  u_\alpha$ \Iff \, $\beta \notin \cup\{v_\gamma:\gamma < \beta\}$;

\item[(e)]   if $\alpha \notin \cup\{v_\gamma:\gamma < \alpha\}$, \then \,:
\begin{enumerate}
\item[$\bullet_1$]  $\gamma \in v_\alpha$ \Iff \, 
($ \gamma > \alpha $ and) 
$G _\gamma $ is embeddable into $G'_\alpha$ over $H$

\item[$\bullet_2$]   if $\gamma \in \alpha_\ast \backslash
 (\alpha +1) \backslash ( \cup\{v_\beta:\beta \le \alpha\ )$ 
then $G_\gamma$ is not embeddable over
  $H$ into any $G'$ satisfying $G'_\alpha \subseteq G' \in \mathbf K_{}$;
\end{enumerate}

\item[(f)]   if $\alpha = \beta +1$ and $\beta \notin
  u_\alpha$ then $v_\beta = \emptyset$.
\end{enumerate}
\end{enumerate}

\underline{Why  this suffices?}

Because if we let $u_{\alpha_\ast} = \alpha_\ast \backslash
(\cup\{v_\gamma:\gamma < \alpha_\ast\})$, 
then $\langle G'_\alpha:\alpha \in u_{\alpha_\ast}\rangle$ is
as required; but we elaborate.

First, for clause $(B)^+$(c) assume that $ \alpha < \beta$ are from 
$u_{\alpha _\ast} $. As $ \beta \notin v_ \alpha $, by 
 $ (\ast)^2_\alpha (e) \bullet _2$  we know that 
 $ G_ \beta $ is not embeddable into any extension of $ G'_ \alpha $  
    over $ H $;   
 but as $ G _ \beta \subseteq G'_ \beta $  clearly 
 also $ G'_ \beta $ is not embeddable into 
 any extension of $ G'_ \alpha $
 over $ H $. Renaming this means that $ G'_\alpha, G'_\beta $  cannot 
 be amalgamated over $ H $, as promised.
 
 Second   for clause $ (B)^+$(d), let $ \alpha \in u_{\alpha _\ast}$  and we have to prove that the pair $ (G'_ \alpha  , H  )$   is   a universal  $ (\mu, \kappa )$-amalgamation
base  where $ \kappa $ is the cardinality of $ H $.  So assume $ G' \in \mathbf{K} _ \mu $ extends $ G'_ \alpha $;  recall that we are assuming that  $ \langle G_ \alpha \colon \alpha < \alpha _\ast \rangle $ is as in clause (B), hence there are $ \beta < \alpha _\ast$  and an  embedding $ f $ of $ G' $ into $ G_ \beta $ over $ H $; we shall prove that $ \beta = \alpha $ hence (recalling $ G_ \alpha \subseteq G'_ \alpha $)  $ f $ embeds $ G' $ into $ G'_ \alpha $  over $ H $ thus finishing proving (B) $\Rightarrow$ (B)$^{+}$. 
    
    If $ \beta \in u_{\alpha _\ast } \setminus \{ \alpha \} $ then  $ f \upharpoonright G'_ \alpha $ embed  $  G'_ \alpha $ into $ G'_\beta $  over $H,$ a contradiction to $ (B)^+$(c) which we have already proved.
    
    If $ \beta \in \alpha _\ast \setminus u_{\alpha _\ast }$ then for some $ \gamma $ 
    we have $ \beta \in v_ \gamma $  hence $ \gamma  < \beta $ and $ G_ \beta $  
    is embeddable into $ G'_\gamma $ over   
    $ H $; hence $ G' $ is embeddable into $ G'_ \gamma $  
    over $ H $. As in the previous sentence necessarily $ \gamma = \alpha $ 
    and we are done.



\underline{Why can we carry out the induction?}

For $\alpha=0,\alpha$ limit we have nothing to do because $u_\alpha$ is determined by $(\ast)^2_\alpha(b)$ and $(\ast)^2_\alpha(c)$.   For $\alpha =  
\beta +1$, if $\beta \in \bigcup\limits_{\gamma < \beta} v_\gamma$ we have nothing to do, in the remaining case we choose $G'_{\beta,i}
\in \mathbf K^{}_\mu$ by induction on 
$i \in [\alpha,\alpha_\ast]$, increasing continuous
with $i $.  
        For $ i=0 $ let 
   $ G'_{\beta, i} = G_\beta$ 
    and for limit $ i $ let $ G'_{\beta, i } = \cup \{ G'_{\beta, j }: j < i  \} $. 
   Then choose $ G'_{\beta,i+1}$ to make clause (e)  
true. 
That is, first 
if $G'_{\beta,i}$ has an
extension into which $G_i$ is embeddable over $H$, then there is such
an extension of cardinality $\mu $;
 and choose $G'_{\beta  ,i+1}$ as such 
an extension.  

Second, if $ G'_{\beta, i } $ has no extension 
into which $ G_i $ is embeddable over $ H $, then we let $G'_{\beta, i + 1 }= G'_{\beta , i}$. 

Lastly, let $G'_\alpha = G'_{\beta,\alpha_\ast}$ and $u_\alpha = u_\beta
\cup \{\beta\}$ and $v_\alpha = \{i \colon i < \alpha_\ast, \, i \geq \alpha ,i \notin \cup
\{v_\gamma:\gamma < \beta\}$ and $G_i$ is embeddable into $G'_\beta$  
over $H\}$.

\underline{$(B)^+ \Rightarrow (A)$}:

We prove below more: there is something like ``special model", i.e. part
(2) of \ref{b40},  that is $ (B)^+ \Rightarrow (C)^+$.

\underline{$(C)^+ \Rightarrow (C) \Rightarrow (A)$} 

This  
is trivial so we are left with 
    proving  
the following.

 \underline{$ (B)^+ \Rightarrow (C)^+$}: 


Let $\mathbf K^{\slf}_{\mu }$ be the class of $ G$ such that:

\begin{enumerate}
\item[$(\ast)^3_{
    G}$]  
\begin{enumerate}
\item[(a)]   $G \in \mathbf K_{\mu };$

\item[(b)]  if $ H \subseteq G, H \in \mathbf{K} _{< \mu }$  
    then there  are $ {\aleph_0} $-indecomposable 
    $ H_ n \subseteq G$ increasing with $n$ for $ n < \omega $ 
    with union of 
    cardinality $ < \mu $  such that $ H \subseteq 
    \cup \{H_n \colon n < \omega   \} $;  and\footnote{For universal classes the ``and'' can be replaced by ``hence''.}
    there are $ {\aleph_0} $-indecomposable $ G_n \subseteq G$
    for $ n < \omega $ such that $ G_n \in \mathbf{K} _{< \mu } $,
        $ G_n \subseteq G_{n + 1 }$  
    and $ G = \cup \{G_n \colon n < \omega \}; $ 
 
\item[(c)] if $ H \subseteq G$  is $ {\aleph_0} $-indecomposable 
of cardinality $ < \mu $
\then \, the pair $(G,H)$ is an
universal 
$ (\mu , < \mu )$-amalgamation base
(see Definition \ref{a32}(2));

\item[(d)]   if $H \subseteq G$ is
  $\aleph_0$-indecomposable of cardinality $< \mu  ,H \subseteq
  H'  \in \mathbf K^{}_{< \mu},H'$ is
  $\aleph_0$-indecomposable\footnote{The $\aleph_0$-indecomposability
    is not always necessary, but we need it sometimes.}, 
and 
  $ G,H'$ are compatible over $H$ (in $\mathbf
  K^{}_{\le \mu}$), \then \, $H'$
  is embeddable into $G$ over $H$.
\end{enumerate}
\end{enumerate}
%
%
%
%

Now we can finish by proving $ (\ast)_4 + (\ast)_5$ below.

\begin{enumerate}
\item[$(\ast)_4$]  If $G \in \mathbf K^{}_{\le \mu}$ \then \,  
some
 $  H  
 \in \mathbf K^{\slf}_{\bar\mu}$  
 extends $G $ 
\end{enumerate} 

We break the proof to some stages, $ (\ast)_{4.3}$
gives the desired conclusion of $(\ast)_4$.  

\begin{enumerate}  
\item[$(\ast)_{4.0}$]  If $ G \in \mathbf{K} _{\le \mu }$  then  for some $ H, \bar{ H } $ we have:

\begin{enumerate} 
    \item[(a)]  $ G \subseteq H \in \mathbf{K} _{\mu };$
    
    \item[(b)]  $ \bar{ H } = \langle H_n \colon n < \omega \rangle;$
    
    \item[(c)] $ H_ n \subseteq H_ {n + 1 } \subseteq H;$
    
    \item[(d)] $ H = \cup \{ H_n:  n < \omega \};$
    
    \item[(e)] each $ H_n $ is $ {\aleph_0} $-indecomposable of cardinality $\mu$; 
    
    \item[(f)] (not really needed)  when $\bfK = \bfK_{\lf},$ if $ K \subseteq H_ n , |K| \le \partial $ and  $ 2^\partial  \le |H_ n |$  then there is a subgroup $ L $  of $ H_n $ extending $ K $ which is $ \Theta ^\orth_\partial $-indecomposable.
\end{enumerate} 
\end{enumerate} 

[Why? For clauses (a)-(e) by the definition of $\bfK$ being nice. For clause (f) by \ref{b34}(1),(2)].

\begin{enumerate} 
    \item[$(\ast)_{4.1}$] if $ N_ 1 \in \mathbf{K} _ {\le \mu} $ \then \, 
    there is $ N_ 2 $ such that 
    \begin{enumerate} 
        \item[(a)] $ N_2 \in \mathbf{K} _\mu;$
        
        \item[(b)] $ N_ 1 \subseteq N_2;$
        
        \item[(c)] if $ H \in \IDC_{< \mu}(N_1)$  \then \, 
        $ (N_2, H) $ is a 
        universal $(\mu, {<} \, \mu )$-amalgamation base.
    \end{enumerate} 
\end{enumerate}

Why? by \ref{b37} it is enough to deal with one such $ H$,  which is O.K. by clause (d) of Def \ref{b28}, recalling $\lqq $universal $(\mu, < \mu )$-amalgamation base"  by (B)$^+$  which we are assuming.

\begin{enumerate} 
    \item[$(\ast)_{4.2}$] like $ (\ast)_{4.1}$ but in clause (c) 
    is replaced by: 
    \begin{enumerate} 
        \item[(c)']  if $ H_1 \in \IDC_{< \mu}(N_1)$   and 
        $ H_1 \subseteq H_2  \in \mathbf{K} _{< \mu }$  (and, we may add, $ H_2 $ is $ {\aleph_0} $-indecomposable) 
        then  \underline{either} $ N_2, H_1 $ are incompatible over $ H_ 1 $ in 
        $ \mathbf{K} _{\le \mu }$  \underline{or} $ H_2 $ is embeddable into $ N_2 $ over $ H_1.$
    \end{enumerate} 
\end{enumerate} 

[Why? Again it is enough to deal with one
pair $ (H_1, H_2)$]  which is done by hand.]

\begin{enumerate} 
    \item[$(\ast)_{4.3}$] If $ N_1 \in \mathbf{K} _{\le \mu }$
    then there is $ N_ 2$  
    such that 
    \begin{enumerate} 
    \item[(a)] $ N_2 \in \mathbf{K} _\mu;$
    
    \item[(b)] $ N_1 \subseteq N_2;$
    
    \item[(c)]  if $ H \in \IDC_{{<} \, \mu}(N_2)$  \then \,  $ (N_2, H) $ is a 
    universal $(\mu, < \mu )$-amalgamation base;
    
    \item[(d)]   if $ H_1 \in \IDC_{{<} \, \mu}(N_2)$   and $ H_1 \subseteq H_2  \in \mathbf{K} _{< \mu }$  (and, we may add, $ H_2 $ is $ {\aleph_0} $-indecomposable) 
    then  \underline{either} $ N_2, H_1 $ are incompatible over $ H_ 1 $ in  $ \mathbf{K} _{\le \mu }$  \underline{or} $ H_2 $ is embeddable into $ N_2 $ over $ H_1.$
    \end{enumerate} 
\end{enumerate} 

[Why? We choose $ L_ \varepsilon \in \mathbf{K} _\mu $ by  induction on $ \varepsilon \le \cf((\mu )$, such that 

\begin{enumerate} 
    \item[(a)]   $ L_ \alpha \in \mathbf{K} _ \mu;$
    
    \item[(b)] $ \langle L_ \beta \colon \beta \le \alpha \rangle$  is increasing continuous;
    
    \item[(c)] $ G_1 \subseteq L_0;$
    
    \item[(d)] if $ \alpha = 3 \beta   + 1  $ then 
    $ L_ \alpha $ relate to $ L_{3 \beta }$ as $ N_2 $ relate to $ N_ 1 $ in $ (\ast)_{4.0}$;
    
    \item[(e)]    if $ \alpha = 3 \beta + 2 $ then 
      $ L_ \alpha $ relate to $ L_{3 \beta  + 1  }$
      as $ N_2 $ relate to $ N_ 1 $ in $(\ast)_{4.1}$;
      
    \item[(f)] if $ \alpha = 3 \beta + 3 $ then $ L_ \alpha $ relate to $ L_{3 \beta   + 2 }$
    as $ N_2 $ relate to $ N_ 1 $ in $ (\ast)_{4.2}.$
\end{enumerate} 

There is no problem to carry the induction and 
note that: if $ N \subseteq L_ {\cf(\mu )}$ is $ \cf(\mu)$-indecomposable the for some $ \varepsilon < \cf(\mu )$ we have $ N \subseteq L_ \varepsilon $. Now  then $ N_2 = L_{\cf(\mu )}$ is as required in $ (\ast)_{4.3}$ hence in $ (\ast)_4$.

\begin{enumerate} 
    \item[$(\ast)_5$]
    \begin{enumerate}
        \item[(a)] if $G_1, G_2 \in \mathbf
          K^{\slf}_{\mu }$ then
        $ G_1, G_2$ are isomorphic;

        \item[(b)]   if $ G_1, G_2 \in \mathbf K^{\slf}_{\mu }, H \in \mathbf K^{}$ is $\aleph_0$-indecomposable and $f_\ell$ embeds $H$ into $G_{\ell} $, for $\ell=1,2$, and this diagram 
        can be completed, (i.e. there are $G \in \mathbf{K} _\mu $ and embedding
        $g_\ell \colon  G_\ell \rightarrow G_\ast$ such that $g_1 \circ f_1 = g_2 \circ f_2$) \then \, there is $h$ such that:

        \begin{enumerate} 
            \item[$(\alpha)$]  $h$ is an isomorphism from $ G_1$ onto $G_2$;
            
            \item[$(\beta)$]  $h \circ f_1 = f_2$;
        \end{enumerate} 
        %
        %
    \end{enumerate}
\end{enumerate}

Why? Clause (a) follows from clause (b) using as $ H $  the trivial group. 
For clause (b), let   
$ {\mathscr F} = {\mathscr F} [G_1, G_2]$ 
be the set of $ f $ such that:

\begin{enumerate} 
    \item[(a)] $ f $ is an isomorphism from  $ G_{1,f} \in \IDC _{{<} \, \mu}(G_1)$ onto $ G_{2,f} \in \IDC_{{<} \, \mu}(G_2);$
    
    \item[(b)] $ G_1, G_2 $ are $ f $-compatible in  $ \mathbf{K} _\mu $
    which means that there is $ G \in \mathbf{K} _ \mu $ and embeddings
    $ g_ {\ell} $ of $ G_ {{\ell} } $ into $ G $ for $ {\ell} = 1,2 $ 
    such that $ g_ 2 \circ f = g_1 \rest G_{1, f } $.
\end{enumerate} 

First  $ {\mathscr F} $ is non-empty (the function $ f$  with domain $ f_1(H) $   and range   $ f_2(H)$  will do). Second use the hence and forth argument; here we use $ \cf(\mu ) = {\aleph_0} $.
\end{PROOF}

\begin{remark}\label{b50}  
    1) Can we prove for strong limit singular $ \mu $ of uncountable cofinality $ \kappa $ a parallel result? Well, we have to consider the following game:
    
    \begin{enumerate}
        \item[$(\ast)$] the game is defined by:
        \begin{enumerate} 
            \item[(a)] a play last $ \theta $ moves,
            
            \item[(b)] in the $ \varepsilon$-th move, first  Player I choose $ M_ \varepsilon \in \mathbf{K} _{ <  \mu} $  and then  player II choose $  N_ \varepsilon  \in \mathbf{K} _{< \mu },$  
            \item[(c)] $M_ \varepsilon \in \mathbf{K} _{< \mu} $ and if $ \varepsilon $ is non-limit then $ M_ \varepsilon $ is $ \cf(\mu ) $-indecomposable,
            
            \item[(d)] $ \langle M_ \zeta \colon \zeta \le \varepsilon \rangle $  is increasing continuous,
            
            \item[(e)] $ M_ \varepsilon \subseteq N_ \varepsilon \subseteq M_{\varepsilon + 1},$
            
            \item[(f)] in the end of the play, the player II wins \Iff \, for every limit ordinal $ \varepsilon < \cf(\mu ) , M_ \varepsilon$  is an amalgamation base inside $ \mathbf{K} _{< \mu }.$
        \end{enumerate} 
    \end{enumerate} 
    
    Now, if player II does not lose then we can imitate the proof above;
    this should be clear.  
    Does the existence of a universal member of $ \mathbf{K} _ \mu $  implies
    this?  we hope to return to this elsewhere.   
    
    See below. 
    
    2) The proof works for any a.e.c. $\mathfrak{k} $ with $\LST_ \mathfrak{k} < \mu $.  But we may wonder: can we  weaken the demand on $ \mathfrak{k} $. Actually, we can: there is no need of smoothness (that is: if $ \langle  M_ \alpha \colon \alpha \le \delta \rangle $  is $ \le _ \mathfrak{k} $-increasing  then  $ \cup \{ M_ \alpha \colon \alpha < \delta \} \le_ \mathfrak{k}  M_ \delta  $). Moreover, while we need the existence of an upper bound for any $ \le _ \mathfrak{k} $-increasing sequence, also we demand the union being such upper bound, only for  the cofinality $ \cf(\mu ) $.

    3)   
    We may add a version fixing $ \bar{ \lambda } $
    
\end{remark} 

We may add (after the journal version): 

\begin{definition}\label{b53}
    We say $\bfK$ is \emph{$\mu$-very nice} when: 

    \begin{enumerate}
        \item[(A)] $\bfK$ is $\mu$-nice (see Definition \ref{b28}),

        \item[(B)] if $\cf(\mu) > \aleph_{0},$ \underline{then} for a club of $\chi < \mu$ we have that $\bfK$ is $\chi$-nice. 
    \end{enumerate}
\end{definition}

\begin{claim}\label{b56}
    If $\bfK$ is $\mu$-very nice then the parallel of Theorem \ref{b40} holds. 
\end{claim}

\newpage

\section{Universal in $\beth_\omega$}\label{3} 

In \S2 we have characterized when there are special models in $\mathbf K$ of cardinality,
e.g., $\beth_\omega$.  We try to analyze a related combinatorial problem. Our intention is to first investigate $\gk_{\fnq}$ (the class structures consisting of a set and  a directed family of equivalence relations on it, each with a finite bound on the size of equivalence classes).  So $\gk_{\fnq}$ is similar to $\mathbf K_{\lf}$ but seems easier to analyze.  We consider some partial orders on $\gk = \gk_{\fnq}$.

First, under the substructure order, $\le_1 \,  = \, \subseteq$, this class fails amalgamation. Second, we have other orders: $\leq_{3}$ demanding  a Tarski-Vaught condition (see below $\TV$).  However, using $ \le _ 3$  where we have a similar demand  for countably many points, finitely many equivalence relations, we
have amalgamation.


This is naturally  connected to locally finite groups, see \ref{b70}, \ref{b73}. 



\begin{definition}\label{b55}
    Let $\mathbf K = \mathbf K_{\fnq}$ be the class of structures $M$ such that (the vocabulary is defined implicitly and is $\tau_{\mathbf K}$, i.e. depends just on $\mathbf K$):
    
    \begin{enumerate}
        \item[(a)]  $P^M,Q^M$ is a partition of $M,P^M$ non-empty;
        
        \item[(b)]  $E^M \subseteq P^M \times P^M \times Q^M$ (is a three-place relation) and we write $a E^M_c b$ for $(a,b,c) \in E^M$;
        
        \item[(c)]  for $c \in Q^M,E^M_c$ is an equivalence relation on $P^M$ with $\sup\{|a/E^M_c|:a \in P^M\}$ finite (see more later);
        
        \item[(d)]  $Q^M_{n,k} \subseteq (Q^M)^n$ for $n,k \ge 1$
        
        \item[(e)]  if $\bar c = \langle c_\ell:\ell < n\rangle \in {}^n(Q^M)$ we
        let $E^M_{\bar c}$ be the closure of $\bigcup\limits_{\ell} E_\ell$ to an equivalence relation; but $E_{\bar{c}}^{M}$ is not part of the vocabulary;
        
        \item[(f)] ${}^n(Q^M) = \bigcup\limits_{k \ge 1} Q^M_{n,k}$;
        
        \item[(g)]  if $\bar c \in Q^M_{n,k}$ then $k \ge |a/E^M_{\bar c}|$ for every $a \in P^M$.
    \end{enumerate}
\end{definition}

\begin{definition}\label{b58}
    We define some partial order on $\mathbf K$.
    
    1) $\le_1 \, = \, \le^1_{\mathbf K} \, = \, \le^1_{\fnq}$ is being a sub-model.
    
    2) $\le_3 = \le^3_{\mathbf K} = \le^3_{\fnq}$ is the following: $M \le_3
    N$ iff:
    
    \begin{enumerate}
        \item[(a)]  $M,N \in \mathbf K,$
        
        \item[(b)]  $M \subseteq N,$
        
        \item[(c)]  if $A \subseteq N$ is countable and $A \cap Q^N$ is finite, \then \, there is an embedding of $N \rest A$ into $M$ over
          $A \cap M$ or just a one-to-one homomorphism.
    \end{enumerate}
    
    3) $\le_2 \, = \, \le^2_{\mathbf K} \, = \, \le^2_{\fnq}$ is defined like $\le_3$
    but in clause (c), $A$ is finite.
\end{definition}
    
\begin{claim}\label{b61}
    1) $\mathbf K$ is a universal class, so $(\mathbf K,\subseteq)$ is an a.e.c.
    
    2) $\le^3_{\mathbf K},\le^2_{\mathbf K},\le^1_{\mathbf K}$ are partial orders
      on $\mathbf K$.
    
    3) $(\mathbf K,\le^2_{\mathbf K})$ is an a.e.c.
    
    4) $(\mathbf K,\le^3_{\mathbf K})$ has disjoint amalgamation.
    
    5) If $M \leq_{2} N, $ $c \in Q^{M}$ and $a \in P^{M},$ then $a / E_{c}^{N}$ is included in $M.$ 

    6) For every $n, k$ there is an existential first order sentence defining, for $M \in K,$ the set $\{ \bar{a} \in {}^{n+2} M \colon a_{n}, a_{n +1} \text{ are } E_{\bar{a}}^{M}\text{-equivalent} \}.$ 
\end{claim}

\begin{PROOF}{\ref{b61}}
    1),2),3) Easy.
    
    4) By \ref{b64} below.
\end{PROOF}

\begin{claim}\label{b64}
    If $M_0 \le^1_{\mathbf K} M_1,M_0 \le^3_{\mathbf K} M_2$ and $M_1 \cap M_2
    = M_0$, \then \, $M = M_1 + M_2$, the disjoint sum of $M_1,M_2$
    belongs to $\mathbf K$ and extends $M_\ell$ for $\ell=0,1,2$ and even
    $M_1 \le^3_{\fnq} M$ and $M_0 \le^2_{\mathbf K} M_1 \Rightarrow M_2
    \le^2_{\mathbf K} M$ \when:
    
    \begin{enumerate}
        \item[$(\ast)$]  $M = M_1 +_{M_0} M_2$ means  $M$ is defined by:  
        
        \begin{enumerate}
            \item[(a)]  $|M| = |M_1| \cup |M_2|$;
            
            \item[(b)]  $P^M = P^{M_1} \cup P^{M_2}$;
            
            \item[(c)]  $Q = Q^{M_1} \cup Q^{M_2}$;
            
            \item[(d)]  we define $E^M$ by defining $E^M_c$ for $c \in Q^M$ by
              cases:
            \begin{enumerate}
                \item[$(\alpha)$]  if $c \in Q^{M_0}$ then $E^M_c$ is the closure of
                  $E^{M_1}_\ell \cup E^{M_2}_\ell$ to an equivalence relation;
                
                \item[$(\beta)$]  if $c \in Q^{M_\ell} \backslash Q^{M_0}$ and $\ell
                  \in \{1,2\}$ then $E^M_c$ is defined by
                \begin{itemize}
                    \item  $a E^M_c b$ \Iff \, $a = b \in P^{M_{3-\ell}} \backslash M_0$ or
                    $a E^{M_\ell}_c b$ so $a,b \in P^{M_\ell}$;
                \end{itemize}
            \end{enumerate}
            
            \item[(e)]  $Q^M_{n,k} $ is the union of $ Q^{M_1}_{n,k} ,  Q^{M_2}_{n,k} $ and the set of 
            $ \bar  {c} $ satisfying 
            \begin{enumerate} 
                \item[$(\alpha )$] $ \bar c \in {}^n(Q^M),$
                
                \item[$(\beta )$] $ \bar{ c } \notin ({}^n(Q^{M_1})) \cup
                  {}^n(Q^{M_2})\}$,
                  
                \item[$(\gamma)$] $E_{\bar{c}}^{M}$ which is now well defined, has no equivalence classes with more than $k$ members, that is, for some finite set $A$ and pairwise distinct $a_{0}, \dots, a_{k} \in A$ which are members of $a / E_{\vec{c}}^{M}$ and the  closure of $\bigcup \{ E_{c_{i}}^{M} \, {\rest} \, A \colon i < \lg(\bar{c}) \}$ to an equivalence relation satisfies $a_{i} E' a$ for $i \leq k$.   
                 
                
                    
            \end{enumerate} 
        \end{enumerate}
    \end{enumerate}
\end{claim}

\begin{PROOF}{\ref{b64}}
    Clearly $M$ is a well defined structure, extends $M_0,M_1,M_2$ and
    satisfies clauses (a),(b),(c) of Definition \ref{b55}.  There are two
    points to be checked: $a \in P^M,\bar c \in Q^M_{n,k} \Rightarrow
    |a/E^M_{\bar c}| \le k$ and ${}^n(Q^M) = \bigcup\limits_{k \ge 1}
    Q^M_{n,k}$
    
    \begin{enumerate}
        \item[$(\ast)_1$]  if $a \in P^M$ and $\bar c \in Q^M_{n,k}$ \then \,
          $|a/E^M_{\bar c}| \le k$.
    \end{enumerate}
    
    [Why?  If $\bar c \in Q^{M_1}_{n,k} \cup
    Q^{M_2}_{n,k}$ this holds by the definition, so assume $\bar c \in
    Q^{M_\iota}_{n,k}, \iota  \le 2$.   
    If this fails, then there is a finite set $A
    \subseteq M$ such that $\bar c \subseteq A,a \in A$ and  the closure of $\bigcup \{ E_{c_{\ell}}^{M} \colon \ell < \lg(\bar{c}) \}$ to an equivalence relation satisfies: every equivalence class has $\leq k$ members. $N = M
    \rest A$ we have $|a/E^N_{\bar c}| > k$.  By $M_0 \le^1_{\mathbf K}
    M_1,M_0 \le^3_{\mathbf K} M_2$ (really $M_0 \le^2_{\mathbf K} M_2$
    suffice) there is a one-to-one homomorphism $f$ from $A \cap M_2$ into
    $M_0$ over $M_{0} \cap A$.  Let $B' = (A \cup M_1) \cup f(A \cap M_2)$ and $N' = M \rest B$
    and let $g = f \cup \id_{A \cap M_1}$.  So $g$ is a homomorphism from
    $N$ onto $N'$ and $g(a)/E^{N'}_{g(\bar c)}$ has $>k$ members, which
    implies $g'(a)/E^{M_1}_{g'(\bar c)}$ has $>k$ members.  Also $g(\bar
    c) \in Q^{M_1}_{n,k}$.  (Why?  If $\iota=1$ trivially, if $\iota=2$ by
    the choice of $f$), contradiction to $M_{1} \in \mathbf K$.]
    
    \begin{enumerate}
        \item[$(\ast)_2$]  if $\bar c \in {}^n(Q^M)$ then $\bar c \in
          \bigcup\limits_{k} Q^M_{n,k}$.
    \end{enumerate}
    
    Why?  If $\bar c \in M_1$ or $\bar c \subseteq M_2$, this is obvious
    by the definition of $M$, so assume that they fail.  By the definition
    of the $Q^M_{n,k}$'s we have to prove that $\sup\{a/E^M_{\bar c} \colon a \in
    P^M\}$ is   finite.  Toward contradiction assume this fails for each
    $k \ge 1$ hence there is $a_k \in P^M$ such that $a_k/E^M_{\bar c}$ has $\ge
    k$ elements hence there is a finite $A_k \subseteq M$ such that 
    $a_k/E^{M \rest A_k}_{\bar c}$ has $\ge k$ elements.  Let $A =
    \bigcup\limits_{k \ge 1} A_k$, so $A_k$ is a countable subset of $M$
    and we continue as in the proof of $(\ast)_1$.
    
    Additional points (not really used) are proved like $(\ast)_2$:
    
    \begin{enumerate}
    \item[$(\ast)_3$]  $M_1 \le^3_{\mathbf K} M$;
    
    \item[$(\ast)_4$]  $M_0 \le^2_{\mathbf K} M_1 \Rightarrow M_2 \le^2_{\mathbf K} M$;
    
    \item[$(\ast)_5$]  $M_1 +_{M_0} M_2$ is equal to $M_2 +_{M_0} M_1$. \qedhere
    \end{enumerate}
\end{PROOF}

\begin{claim}\label{b67}
    1) If $\lambda = \lambda^{< \mu}$ and $M \in \mathbf K$ has cardinality
    $\le \lambda$ \then \, there is $N$ such that:
    
    \begin{enumerate}
        \item[(a)]  $N \in \mathbf K_\lambda$ extend $M$;
        
        \item[(b)]  if $N_0 \le^3_{\mathbf K} N_1$ and $N_0$ has cardinality $<
          \mu$ and $f_0$ $\leq_{2}$-embeds $N_0$ into $N$, \then \, there is an embedding
          $f_1$ of $N_1$ into $N$ extending $f_0$.
    \end{enumerate}

    2) For every $ M \in \mathbf{K} $ we can define an equivalence relation
    $ E = E_ \mathbf{K} $ on the class $ \{N \in \mathbf{K} \colon M \le _2 N  \} $ 
    with  $ \le 2^{{\|M\|}^ {\aleph_0}} $-equivalence  classes such that;
    if $ N_1, N_2 $  are $ E $ -equivalence then they can be amalgamated
    over $ M $ (in $ (\mathbf{K}, \le _ 2) $.

    3) If $ \mu $ is strong limit 
    then $ (\mathbf{K}, \le _2)$  is $ \mu $-nice.
\end{claim}

What is the connection to $ \mathbf{K} _{\lf}$?  the following   
explain (see \cite{KeWe73})

\begin{definition} \label{b70}
    1) For a group $ G \in \mathbf{K} _{\lf} $ we define $ M= \fnq_G \in \mathbf{K} _{\fnq}$  as follows:
    
    \begin{enumerate}  
        \item[(a)] $ P^M$ is the set of elements of $ G,$
        
        \item[(b)]  $ Q^M= \{ (c,1) \colon c \in G\} $, a copy of $ G,$
        
        \item[(c)] $ E^M$  is the set of triples $ (a,b,(c,1)$ such that
        $ a,b,c \in G $ and for some $ n,m \in \mathbb{Z} $  we have
        $ G \models c^n a c^m = b $.
    \end{enumerate}

    2) For $ M \in \mathbf{K} $ we define  $ G= \grp _ M$  as the subgroup 
    of $ \sym(P^M)$  consisting of the permutations $ \pi $ of $ P^M$
    such that for some finite sequence $ \bar{ c } $ of elements of $ Q^M $ 
    we have: for every $ x \in P^M$  we have
    $ \pi(x)  E^M_ {\bar{ c}}  x$. 
\end{definition} 

\begin{discussion}\label{b73}
    The problem is that cases of amalgamation in $ (\mathbf{K} , \le _2)$ cannot
    be lifted to one in $ \mathbf{K} _{\lf},$ that is, in \ref{b64}, for $c \in M_{\ell} \setminus M_{0},$ we can choose $E_{c}^{M} \, {\rest} \, (M_{3-\ell} \setminus M_{0})$ as the equality but the parallel demand for groups fail. 
\end{discussion} 

After publication we add: 

\begin{claim}\label{b76}
    Assume that $\mu$ is a strong limit singular cardinal of cofinality $\aleph_{0}.$ Then $(\bfK_{\fnq}, \leq_{2/3})$ has a universal member in $\mu.$ 
\end{claim}

\begin{PROOF}{\ref{b76}}
    By the general criterion \cite[1.16 = L a34]{Sh:820} + $\JEP,$ but we elaborate. 

    \begin{enumerate}
        \item[$(\ast)_{1}$] fix $\bar{\lambda} = \langle \lambda_{n} \colon n < \omega \rangle,$ $2^{\lambda_{n}} < \lambda_{n + 1} < \mu = \sum \{ \lambda_{k} \colon k < \omega\}$ and $\lambda_{n} = \lambda_{n}^{\aleph_{0}}.$   
    \end{enumerate}

    \begin{enumerate}
        \item[$(\ast)_{2}$] For $\xi \leq \omega,$ let: 

        \begin{enumerate}
            \item[(a)] $\bfK_{\xi} = K_{\bar{\lambda}, \xi}$ is the class of $\bar{M}$ such that $\bar{M} = \langle M_{n} \colon n < \xi \rangle,$ $M_{n} \in \bfK_{\lambda_{n}}^{\fnq}$ is $\leq_{3}$-increasing with $n,$

            \item[(b)] $\bfK_{\xi} = \bfK_{\bar{\lambda}, \xi} = \{ \bar{M} \in \bfK_{\xi} \colon$the universe of $M_{n}$ is $\lambda_{n}$ for any $n < \xi\}.$ 
        \end{enumerate}
    \end{enumerate}

    \begin{enumerate}
        \item[$(\ast)_{3}$] if $M \in \bfK_{M}^{\fnq}$ then there is $\bar{M} \in \bfK_{\bar{\lambda}, \omega}$ whose union is $M$.  
    \end{enumerate}

    \begin{enumerate}
        \item[$(\ast)_{4}$] We can find $\bar{N}$ such that: 

        \begin{enumerate}
            \item[(a)] $\bar{N} = \left \langle N_{\alpha, \eta} \colon \alpha < 2^{\lambda_{0}}, \, \eta \in \prod_{\ell < n} 2^{\lambda_{\ell + 1}} \text{ for some } n < \omega \right \rangle,$   

            \item[(b)] $\langle N_{\alpha, \langle \, \rangle} \colon \alpha < 2^{\lambda_{0}} \rangle $ list the $M \in \bfK_{\fnq}$  with universe $\lambda_{0},$

            \item[(c)] for $\alpha < 2^{\lambda_{0}},$ $\eta \in \prod_{\ell < n}2^{\lambda_{\ell + 1}},$ the sequence $\langle N_{\alpha, \eta^{\smallfrown} \langle \beta \rangle} \colon \beta < 2^{\lambda_{ n+1}}$ list the $M \in \bfK_{\fnq}$ with universe $\lambda_{n+1}$ which $<_{3}$-extend $N_{\alpha +1}.$  
        \end{enumerate}
    \end{enumerate}

    \begin{enumerate}
        \item[$(\ast)_{5}$] We can find $N_{\ast}, \bar{h}$ such that: 

        \begin{enumerate}
            \item[(a)] $N_{\ast} \in \bfK_{\fnq}$ has cardinality $\mu,$

            \item[(b)] $\bar{h} = \left \langle h_{\alpha, \eta} \colon \alpha < 2^{\lambda_{0}}, \, \eta \in \prod_{\ell < n} 2^{\lambda_{\ell +1}} \text{ for some } n < \omega 
            \right \rangle,$

            \item[(c)] $h_{\alpha, \eta}$ embeds $N_{\alpha, \eta}$ into $N_{\ast},$

            \item[(d)] if $\nu \lhd \eta \in \prod_{\ell < n} 2^{\lambda_{\ell +1}} $ and $ \alpha < 2^{\lambda_{0}},$ then $h_{\alpha, \eta} \subseteq h_{\alpha, \nu}.$
        \end{enumerate}
    \end{enumerate}

    \begin{enumerate}
        \item[$(\ast)_{6}$] $N_{\ast}$ is a universal member of $\bfK_{\fnq}$ in $\mu.$ 
    \end{enumerate}

    [Why? By $(\ast)_{3}$ and $(\ast)_{5}$.]
\end{PROOF}

\newpage


\bibliographystyle{amsalpha}
\bibliography{shlhetal}

\end{document}